\documentclass{article}

\usepackage{arxiv}

\usepackage[utf8]{inputenc} 
\usepackage[T1]{fontenc}    
\usepackage{hyperref}       
\usepackage{url}            
\usepackage{booktabs}       
\usepackage{amsfonts}       
\usepackage{nicefrac}       
\usepackage{microtype}      
\usepackage{lipsum}

\title{Non-linear extension of interval arithmetic and exact resolution of interval equations over square regions}

\author{
  Giovanny A.~Fuentes S \\
  Universidade Federal Fluminesi \\
  Niteroi,RJ\\Brasil\\
  \textit{giovannyfuentes@id.uff.br} \\
}

\begin{document}
\maketitle

\begin{abstract}
The interval numbers is the set of compact intervals of $\mathbb{R}$ with addition and multiplication operation, which are very useful for solving calculations where there are intervals of error or uncertainty, however, it lacks an algebraic structure with an inverse element, both additive and multiplicative This fundamental disadvantage results in overestimation of solutions in an interval equation or also overestimation of the image of a function over square regions. In this article we will present an original solution, through a morphism that preserves both the addiction and the multiplication between the space of the interval numbers to the space of square diagonal matrices.
\end{abstract}


\tableofcontents
\section{Introduction}
In this article we propose a general method to solve equations with interval variables, that is to say where the unknown of the problem  is a compact interval, through a $\varphi$ application that takes intervals and takes them to 2 \ times 2 square matrices. In the first part of this article we will give a brief introduction to the Artificial Interval. In the second part we will state and demonstrate \textit{the fundamental theorem of the interval functions}, which gives us the necessary conditions so as not to have lost points, when computing the image of a square region (finite Cartesian producer of compact intervals) through the $\varphi$ application. And finally in the third and last part of this article, through \textit{the fundamental theorem of the interval functions}, we will state and demonstrate the fundamental theorem of the interval equations, that will give us necessary and sufficient conditions to solve the interval equations.
\subsection{Basic Terms and Concepts the interval arithmetic }
Recall that the closed interval demoted by $[a,b]$ is the real numbers given by 
\begin{align}
[a,b]=\{x\in\mathbb{R}, a\leq x\leq b\}
\end{align}
We say that an interval is Degenerate if $a=b$. Such an interval contains a single real number $a$. By convention, we agree to identify a degenerate interval $[a,a]$. In this sense, we may write such equation as
\begin{align}
0=[0,0]
\end{align}
We will denote by $K_c(\mathbb{R})$ the set of compact intervals real.\\\\
We are about to define the basic arithmetic operations between intervals. The key point in these definition is that computing  with set. For example when we add two interval, the resulting interval is set containing the sums of all pair of number, one form each of the initial sets. By definition then, the sum of two intervals $X$ and $Y$ is the set
\begin{align}
X+Y=\{ x+y; x\in X, y\in Y\}
\end{align}
The difference of two intervals $X$ and $Y$ is the set
\begin{align}
X-Y=\{x-y; x\in X. y\in Y\}
\end{align} 
The product of $X$ and $Y$ is given by 
\begin{align}
XY=\{xy;x\in X, y\in Y\}.
\end{align}
Finally, the quotient $X/Y$ with $0\not\in Y$ is defined as
\begin{align}
X/Y=\{x/y;x\in X, y\in Y   \}
\end{align}
\subsection{Endpoint formulas for the arithmetic Operations}
\textbf{Addition} Let us that an operational way to add intervals. Since $x\in X=[x_1,x_2]$ means that $x_1\leq x\leq x_2$ and $y\in Y=[y_1,y_2]$ means that  $y_1\leq y\leq y_2$ , we see by addition of inequalities that the numerical sums $x+y\in X+Y$ must satisfy $x_1+y_1\leq x+y\leq x_2+y_2$. Hence, the formula $X+Y=[x_1+y_1,x_2+y_2]$
\begin{example} Let $X=[0,2]$ and $Y=[-1,1]$. Then $X+Y=[0-1,1+2]=[-1,3]$
\end{example}
\textbf{Subtraction} Let $X=[x_1,x_2]$ and $Y=[y_1,y_2]$.  We add the inequalities
\begin{align}
x_1\leq x\leq x_2\text{ and } -y_2\leq -y\leq -y_1
\end{align}
to get $x_1-y_2\leq x-y\leq x_2-y_1$. It follows that $X-Y=[x_1-y_2,x_2-y_1]$. Note that $X-Y=X+(-Y)$ where $-Y=[-y_2,-y_1]$\\
\begin{example} Let $X=[-1,0]$ and $Y=[1,2]$. Then $X-Y=[-1-2,0-1]=[-3,-1]$
\end{example}
\textbf{Multiplication} In terms of endpoint, the product $XY$ of two intervals $X$ and $Y$ is given by 
\begin{align}
XY=[\min S,\max S], \text{ where } S=\{ x_1y_1,x_1y_2,x_2y_1,x_2y_2\}
\end{align}
 \begin{example} Let $X=[-1,0]$ and $Y=[1,2]$. Then $S=\{-1,-2,0\}$ and  $XY=[-2,0]$.
 \end{example}
The multiplication of intervals is given in terms of the minimum and maximum of  four products of endpoint, this can be broken into nine spacial cases. Let $X=[x_1,x_2]$ and $Y=[y_1,y_2]$ then $XY=Z=[z_1,z_2]$, then
\begin{center}
\begin{tabular}{|c|c|c|}
\hline 
case & $z_1$ & $z_2$ \\ 
\hline 
$0\leq x_1,y_1$ & $x_1y_1$ & $x_2y_2$ \\ 
\hline 
$x_1<0<x_2$ and $0\leq y_1$ & $x_1y_2$ & $x_2y_2$ \\ 
\hline 
$x_2\leq 0$ and $0\leq y_1$ & $x_1y_2$ & $x_2y_1$ \\ 
\hline 
$0\leq x_1$ and $y_1<0< y_2$ & $x_2y_1$ & $x_2y_2$ \\ 
\hline 
$x_2\leq 0$ and $y_1<0< y_2$ & $x_1y_2$ & $x_1y_1$ \\ 
\hline 
$0\leq x_1$ and $y_2\leq 0$ & $x_2y_1$ & $x_1y_2$ \\ 
\hline 
$x_1< 0<x_2$ and $y_2\leq 0$ & $x_2y_1$ & $x_1y_1$ \\ 
\hline 
$x_2\leq 0$ and $y_2\leq 0$ & $x_2y_2$ & $x_1y_2$ \\ 
\hline 
$x_1<0<x_2$ and $y_1<0<y_2$ & $\min\{x_1y_2,x_2y_1\}$ & $\max\{x_1y_1,x_2y_2\}$ \\ 
\hline 
\end{tabular}
\end{center}
\vspace{0.5cm}
\textbf{Division} As with real number, division can accomplished via multiplication by the reciprocal of the second operand. That is, we can implement equation using
\begin{align}
X/Y=X\left( \frac{1}{Y}\right) ,
\end{align}
where 
\begin{align}
\frac{1}{Y}=\left\lbrace y;\frac{1}{y}\in Y\right\rbrace .
\end{align}
Again, this assume $0\not\in Y$. For more information about interval numbers, see \cite{moore}\\\\

\section{Interval Function Several Variable  and his Fundamental Theorem}
The main result of this section is the fundamental theorem of the interval functions, which allows us to calculate the value of a function on a square region, that is, on a finite Cartesian product of compact intervals, this can be interpreted as determining the image of a Interval function where each variable detects a unique value under the same symbol.
Consider $\mathcal{M}_{2x2}(\mathbb{R}^n)$ as the space of the square matrices two by two over $\mathbb{R}^n$. Let $\Vert.\Vert_n$ a norm in  $\mathbb{R}^n$ and $\Vert.\Vert_s$ the norm of sum in $\mathcal{M}_{2x2}(\mathbb{R})$   then we define the following norm in $\mathcal{M}_{2x2}(\mathbb{R}^n)$ as 
\begin{align}
\left\Vert \left( \begin{array}{cc}
 (a_1,\ldots,a_n) &  (b_1,\ldots,b_n)  \\ 
 (c_1,\ldots,c_n)  &  (d_1,\ldots,d_n) 
\end{array} \right)\right\Vert= 
\Vert(a_1,\ldots,a_n)\Vert_n+\Vert(b_1,\ldots,b_n)\Vert_n +
\Vert(c_1,\ldots,c_n)\Vert_n +\Vert(d_1,\ldots,d_n)\Vert_n
\end{align}
We have this application, it satisfies the norms of norm, then the space $\mathbb{R}^n$ is a normed space.
  Let $\Vert .\Vert$ be a matrix norm in $\mathcal{M}_{2x2}(\mathbb{R}^n)$. In this section we will consider $\mathcal{D}_2(\mathbb{R}^n)$ the space of the square diagonal matrices two by two over $\mathbb{R}^n$  with the following norm $\Vert.\Vert_*=\frac{1}{2}\Vert .\Vert$. Obviously $\Vert.\Vert_*$ is a norm in $\mathcal{D}_2(\mathbb{R}^n)$.
\begin{proposition}
$\mathcal{M}_{2x2}(\mathbb{R}^n)$ with a norm defined is complete.
\end{proposition}
\textbf{Proof} Let $\{A_j\}_j$ an Cauchy sequence in $\mathcal{M}_{2x2}(\mathbb{R}^n)$ and let $\varepsilon>0$, then exist a $N>0$ such that for $p,m>N$ we have $\Vert A_p-A_m\Vert<\varepsilon$, then 
\begin{align}
&\left\Vert \left( \begin{array}{cc}
 (a_{p1},\ldots,a_{pn})-(a_{m1},\ldots,a_{mn}) &  (b_{p1},\ldots,b_{pn})-(b_{m1},\ldots,b_{mn})   \\ 
 (c_{p1},\ldots,c_{pn})-(c_{m1},\ldots,c_{mn})   &  (d_{1p},\ldots,d_{pn})-(d_{m1},\ldots,d_{mn})  
\end{array} \right)\right\Vert\\
&=\Vert(a_{p1},\ldots,a_{pn})-(a_{m1},\ldots,a_{mn})\Vert_n+\Vert(b_{p1},\ldots,b_{pn})-(b_{m}1,\ldots,b_{mn})\Vert_n\\ +
&\Vert(c_{p1},\ldots,c_{pn})-(c_{m1},\ldots,c_{mn})\Vert_n +\Vert(d_{1n},\ldots,d_{pn})-(d_{m1},\ldots,d_{mn})\Vert_n<\varepsilon
\end{align}
this implique that $\Vert(x_{p1},\ldots,x_{pn})-(x_{m1},\ldots,x_{mn})\Vert_n<\varepsilon$ for $x=a,b,c,d$, as $\mathbb{R}^n$ is complete, then $(x_{p1},\ldots,x_{pn})$ is convergent. Therefore $\{A_j\}_j$ is convegent. Which proves that $\mathcal{M}_{2x2}(\mathbb{R}^n)$ is complete.
\begin{proposition}\label{matrixconver} Let $\{A_j\}_j$ a sequence in $\mathcal{M}_{2x2}(\mathbb{R}^n)$, If the series $\displaystyle\sum_{j=0}^{\infty}\parallel A_j\parallel$ converges for any norm in, so the series $\displaystyle\sum_{j=0}^{\infty}A_j$ is converges.
\begin{flushright}
$\blacksquare$
\end{flushright}
\end{proposition}
\textbf{Proof} 
Let $\varepsilon>0$, We will see that there is an integer $N>0$  such that $p,q\geq N$, then $\Vert S_p-S_q\Vert<\varepsilon$, where $\displaystyle S_p=\sum_{j=0}^pA_j$ This shows that the succession $S_n$  of partial sums is a Cauchy sequence, As $\mathcal{M}_{2x2}(\mathbb{R}^n)$ it's complete, then $\{S_n\}$ is converges. Indeed:
\begin{align}
\Vert S_p-S_q\Vert=\Vert\sum_{j=p+1}^qA_j\Vert\leq\sum_{j=p+1}^q\Vert A_j\Vert=\left|\sum_{j=0}^p\Vert A_j\Vert-\sum_{j=0}^q\Vert A_j\Vert\right|
\end{align}
Now, if $\displaystyle\sum_{j=0}^\infty\Vert A_j\Vert$ converges, exist $N>0$ such that $p,q\geq N$, then 
\begin{align}
\left|\sum_{j=0}^p\Vert A_j\Vert-\sum_{j=0}^p\Vert A_j\Vert\right|<\varepsilon.
\end{align}
what we wanted to demonstrate.
\begin{flushright}
$\blacksquare$
\end{flushright}
\begin{definition}  Let $f:\mathcal{X}\subset\mathbb{R}^n\to\mathbb{R}$ analytic  function and $a=(a_1,\ldots,a_n)\in int(\mathcal{X})$ and $\varepsilon>0$ such that $B(a,\varepsilon)\subset  int(\mathcal{X})$, and $(\alpha_1,\ldots \alpha_n),(\beta_1,\ldots \beta_n)\in B(a,\varepsilon)$, we defined
\begin{align}
&f\left( \prod_{j=1}^n \left( \begin{array}{cc}
x_j & 0 \\ 
0 & y_j
\end{array} \right)\right)\\ &=\displaystyle\sum_{j=0}^{\infty}\frac{1}{j!}\sum_{j_1+\ldots+j_n=j}\binom{j}{j_1\ldots j_n}\frac{\partial^n f(a_1,\ldots,a_n)}{\partial x_1^{j_1}\ldots\partial x_n^{j_n}}\left(\begin{array}{cc}
a_1-\alpha_1 & 0 \\ 
0 & a_1-\beta_1
\end{array} \right)^{j_1}\ldots\left( \begin{array}{cc}
a_n-\alpha_n & 0 \\ 
0 & a_n-\beta_n
\end{array} \right)^{j_n}
\end{align}
where $\binom{j}{j_1\ldots j_n}$ is an mutinomial coeficent.
\end{definition}
 
 \begin{definition} Let $f:\mathcal{X}\subset\mathbb{R}^n\to\mathbb{R}$ analytic  function and $a=(a_1,\ldots,a_n)\in int(\mathcal{X})$ and $\varepsilon>0$ such that $B(a,\varepsilon)\subset  int(\mathcal{X})$, and $(\alpha_1,\ldots \alpha_n),(\beta_1,\ldots \beta_n)\in B(a,\varepsilon)$, we defined  $\varphi:Kc(\mathbb{R})\to \mathcal{D}_2(\mathbb{R})$   by  $\varphi([a,b])=\left(\begin{array}{cc}
a & 0 \\ 
0 & b
\end{array}\right)$ and $\overline{\varphi}:\mathcal{F}(Kc(\mathbb{R)}^n)\to \mathcal{F}(\mathcal{D}_2(\mathbb{R})^n)$   by 
\begin{equation}
\overline{\varphi} f\left( \prod_{j=1}^n[a_j,b_j]\right):=f\left( \prod_{j=1}^n\varphi [a_j,b_j]\right) =f\left( \prod_{j=1}^n \left( \begin{array}{cc}
a_j & 0 \\ 
0 & b_j
\end{array} \right)\right) 
\end{equation}
\end{definition}
\begin{proposition}\label{analy} Let $f:\mathcal{X}\subset\mathbb{R}^n\to\mathbb{R}$ analytic  function and $a=(a_1,\ldots,a_n)\in\mathbb{R}^n$ such that
$B(x_0,\varepsilon)\subset\mathcal{X}$ for any $\varepsilon>0$ and $\displaystyle\prod_{j=1}^n[\alpha_j,\beta_j]\subset B(x_0,\varepsilon)$ , then
\begin{equation}
f\left( \prod_{j=1}^n\left( \begin{array}{cc}
\alpha_j & 0 \\ 
0 & \beta_j
\end{array} \right)\right) =\left( \begin{array}{cc}
f\left(\displaystyle \prod_{j=1}^n\alpha_j \right)& 0 \\ 
0 & f\left( \displaystyle\prod_{j=1}^n\beta_j  \right)
\end{array} \right) 
\end{equation}
\end{proposition}
\textbf{Proof} Consider any norm $\Vert.\Vert_n$ of $\mathbb{R}^n$. Let $a\in int\mathcal{X}$, as $f$ is analytic in $\mathcal{X}$, then exist $\varepsilon>0$ such that $B(x_0,\varepsilon)\subset\mathcal{X}$ such that for all $x\in B(x_0,\varepsilon)$ there is a series of Taylor that converges for $f(x)$. Let $\displaystyle\prod_{j=1}^n[\alpha_j,\beta_j]\subset B(x_0,\varepsilon)$, consider a matrix
$
\left( \begin{array}{cc}
(a_1,\ldots,a_n) & 0 \\ 
0 & (b_1,\ldots,b_n)
\end{array} \right) 
$
that is matrix representation associated with the interval $\displaystyle\prod_{j=1}^n[\alpha_j,\beta_j]$ and $\left( \begin{array}{cc}
(a_1,\ldots,a_n) & 0 \\ 
0 & (a_1,\ldots,a_n)
\end{array} \right) $ the representation  of  the matrix associated with the interval $\displaystyle\prod_{j=1}^n[a_j,a_j] $. Consider a norm od in space $\mathcal{M}_{2x2}(\mathbb{R}^n)$. Note that:
\begin{align}
\left|\left|
\left( \begin{array}{cc}
(a_1,\ldots,a_n) & 0 \\ 
0 & (a_1,\ldots,a_n)
\end{array} \right) -\left( \begin{array}{cc}
(\alpha_1,\ldots,\alpha_n) & 0 \\ 
0 & (\beta_1,\ldots,\beta_n)
\end{array} \right) 
\right|\right|_{*}\\
=\frac{1}{2}(\Vert(a_1,\ldots,a_n)-(\alpha_1,\ldots,\alpha_n)\Vert_n+\Vert(a_1,\ldots,a_n)-(\beta_1,\ldots,\beta_n)\Vert_n)<\varepsilon.
\end{align}
then by proposition \ref{matrixconver} the series of matrix power
\begin{align}
\displaystyle\sum_{j=0}^{\infty}\frac{1}{j!}\sum_{j_1+\ldots+j_n=j}\binom{j}{j_1\ldots j_n}\frac{\partial^n f(a_1,\ldots,a_n)}{\partial x_1^{j_1}\ldots\partial x_n^{j_n}}\left(\begin{array}{cc}
a_1-\alpha_1 & 0 \\ 
0 & a_1-\beta_1
\end{array} \right)^{j_1}\ldots\left( \begin{array}{cc}
a_n-\alpha_n & 0 \\ 
0 & a_n-\beta_n
\end{array} \right)^{j_n}
\end{align} is convergent.\\\\
Now consider the partial  sum
\begin{align}
&\displaystyle\sum_{j=0}^{m}\frac{1}{j!}\sum_{j_1+\ldots+j_n=j}\binom{j}{j_1\ldots j_n}\frac{\partial^n f(a_1,\ldots,a_n)}{\partial x_1^{j_1}\ldots\partial x_n^{j_n}}\left(\begin{array}{cc}
a_1-\alpha_1 & 0 \\ 
0 & a_1-\beta_1
\end{array} \right)^{j_1}\ldots\left( \begin{array}{cc}
a_n-\alpha_n & 0 \\ 
0 & a_n-\beta_n
\end{array} \right)^{j_n}\\
&=\left( \begin{array}{cc}
\Gamma_m(a_1-\alpha_1)^{j_1}\ldots(a_n-\alpha_n)^{j_n} & 0 \\ 
0 & \Gamma_m(a_1-\beta_1)^{j_1}\ldots(a_n-\beta_n)^{j_n}
\end{array} \right),
\end{align}
where $\Gamma_m=\displaystyle\sum_{j=0}^{m}\frac{1}{j!}\sum_{j_1+\ldots+j_n=j}\binom{j}{j_1\ldots j_n}\frac{\partial^n f(a_1,\ldots,a_n)}{\partial x_1^{j_1}\ldots\partial x_n^{j_n}}$. As $(\alpha_1,\ldots,\alpha_n),(\beta_1,\ldots,\beta_n)\in B(a,\varepsilon)$, then $\Gamma_m(a_1-\alpha_1)^{j_1}\ldots(a_n-\alpha_n)^{j_n}\to f\left(\displaystyle \prod_{j=1}^n\alpha_j \right)$ and $\Gamma_m(a_1-\beta_1)^{j_1}\ldots(a_n-\beta_n)^{j_n}\to f\left(\displaystyle \prod_{j=1}^n\beta_j \right) $ as $m\to\infty$, then 
\begin{equation}
\left( \begin{array}{cc}
\Gamma_m(a_1-\alpha_1)^{j_1}\ldots(a_n-\alpha_n)^{j_n} & 0 \\ 
0 & \Gamma_m(a_1-\beta_1)^{j_1}\ldots(a_n-\beta_n)^{j_n}
\end{array} \right)\to \left( \begin{array}{cc}
f\left(\displaystyle \prod_{j=1}^n\alpha_j \right)& 0 \\ 
0 & f\left(\displaystyle \prod_{j=1}^n\beta_j \right)
\end{array} \right)
\end{equation}
as $m\to\infty$, then
\begin{equation}
f\left( \prod_{j=1}^n\left( \begin{array}{cc}
\alpha_j & 0 \\ 
0 & \beta_j
\end{array} \right)\right) =\left( \begin{array}{cc}
f\left(\displaystyle \prod_{j=1}^n\alpha_j \right)& 0 \\ 
0 & f\left( \displaystyle\prod_{j=1}^n\beta_j  \right)
\end{array} \right) 
\end{equation}
 \begin{flushright}
$ \blacksquare$
 \end{flushright}

\begin{definition} Let $f:\mathcal{X}\subset\mathbb{R}^n\to\mathbb{R}$ differential  function and $\displaystyle\prod_{j=1}^n[a_j,b_j]\subset\mathcal{X}$, we say that $\displaystyle\prod_{j=1}^n[a_j,b_j]$ is free of singularity if the components of the gradient vector are different from zero in all  $\displaystyle\prod_{j=1}^n[a_j,b_j]$, that is to say
\begin{align}
\frac{\partial f(x)}{\partial x_j} \not =0\text{ for all } x\in \displaystyle\prod_{j=1}^n[a_j,b_j]\text{ and } j=1,\ldots,n.
\end{align}
\end{definition}
\begin{definition}\label{switch} Let  $f:\mathcal{X}\subset\mathbb{R}^n\to\mathbb{R}$ analytic function, $\displaystyle\prod_{j=1}^n[a_j,b_j]\subset \mathcal{X}$ free of singularity and $[a_j,b_j]$ an interval of the $j$-variable. Defined the switch $\widehat{[a_j,b_j]}$ of the interval as 
 $[a_j,b_j]$ if $\displaystyle\frac{\partial f}{\partial x_j}>0$ and  $[b_j,a_j]$ if $\displaystyle\frac{\partial f}{\partial x_j}<0$, and denote by $\phi f\left(\displaystyle\prod_{j=1}^n[a_j,b_j]\right)=\overline{\varphi}  f\left(\displaystyle\prod_{j=1}^n\widehat{[a_j,b_j]}\right)$.
\end{definition}
Note that $\phi$ is well defined, since the swapper depends only on the sign of the directional derivative on the square region, and this does not depend on how the function is formulated and also how the square region is free of singularities, then this sign is constant in the square region. And on the other hand, since $f$ is an analytical function, we have that the matrix series is convergent.

\begin{theorem} [Fundamental Theorem of  Interval Functions]\label{teofun} 
Let $f:\mathcal{X}\subset\mathbb{R}^n\to\mathbb{R}$ analytic function and $x_0\in\mathbb{R}^n$ such that
$B(x_0,\varepsilon)\subset\mathcal{X}$ for any $\varepsilon>0$ and $\displaystyle\prod_{j=1}^n[a_j,b_j]\subset B(x_0,\varepsilon)$ free of singularity, then
\begin{align}
  \displaystyle f\left(\prod_{j=1}^n {[a_j,b_j]}\right)=\varphi^{-1}\phi f\left( \prod_{j=1}^n[a_j,b_j]\right)  
\end{align}
  .
\end{theorem}
\textbf{Proof} Let $f:\mathcal{X}\subset\mathbb{R}^n\to\mathbb{R}$ analytic function and $x_0\in\mathbb{R}^n$ such that
$B(x_0,\varepsilon)\subset\mathcal{X}$ for any $\varepsilon>0$ and $\displaystyle\prod_{j=1}^n[a_j,b_j]\subset B(x_0,\varepsilon)$ free of singularity, then 
\begin{align}
\displaystyle\phi f\left( \prod_{j=1}^n\widehat{[a_j,b_j]}\right)=\displaystyle\overline{\varphi}  f\left( \prod_{j=1}^n[x_j,y_j]\right)
\end{align}
where $[x_j,y_j]$   is the result of applying the switch to $[a_j,b_j]$, then by Proposition \ref{analy}  since $f$ is analytic function, we have
\begin{align}
\displaystyle\overline{\varphi}  f\left( \prod_{j=1}^n[x_j,y_j]\right)&=\left( \begin{array}{cc}
f\left(\displaystyle \prod_{j=1}^nx_j \right)& 0 \\ 
0 & f\left( \displaystyle\prod_{j=1}^ny_j  \right)
\end{array} \right) .
\end{align}
Applying $\varphi^{-1}$, we have the following interval,
\begin{align}
\left[  f\left( \prod_{j=1}^nx_j\right),f\left( \prod_{j=1}^ny_j\right)\right] .
\end{align}
Now we will prove that the interval above corresponds to the image of $f$ on $ R=\displaystyle\prod_{j=1}^n[a_j,b_j]$. 
First we observe that both $  f\displaystyle\left( \prod_{j=1}^nx_j\right)$ and $  f\displaystyle\left( \prod_{j=1}^nx_j\right)$ are elements of $f(R)$, since $R$ is connected and closed, we have that $\displaystyle\left[  f\left( \prod_{j=1}^nx_j\right),f\left( \prod_{j=1}^ny_j\right)\right] $ is a subset of $f(R)$.\\\\
On the other hand. As we have that the gradient is different from zero within $R$, we have that the function reaches its extreme points at the border of $R$. First we will prove that the maximum and minimum point are not in the edges of $R$, but are in one of the vertex. Indeed, since the gradient  has non-zero components in all $R$ in particularly in the edge of $R$, we have that the extreme points of $f$ can not be in the edges, since $f$ restricted to the edges, it is a continuous function over a compact interval where its derivative is not null, then its extreme points are the extremes of the interval. Since the argument is valid for all edges, we conclude that the extreme points of $f$ are in the vertex. Let $(z_1,\ldots,z_j,\ldots,z_n)$ be a vertex of $R$. Since $f$ is monotonous on the edges, we have the following inequality in $z_j$ leaving the other variables fixed
\begin{align}
f(z_1,\ldots,x_j,\ldots,z_n)\leq f(z_1,\ldots,z_j,\ldots,z_n)\leq f(z_1,\ldots,y_j,\ldots,z_n).
\end{align}
Taking this inequality inductively on each variable, we have
\begin{align}
f(z_1,\ldots,z_j,\ldots z_n)\geq f(x_1,\ldots,z_j,\ldots z_n) \geq f(x_1,\ldots,x_j,\ldots z_n) \geq f(x_1,\ldots,x_j,\ldots x_n)
\end{align}
and
\begin{align}
f(z_1,\ldots,z_j,\ldots z_n)\leq f(y_1,\ldots,z_j,\ldots z_n)\leq f(y_1,\ldots,y_j,\ldots z_n) \leq f(y_1,\ldots,y_j,\ldots y_n).
\end{align}
Hence
\begin{align}
f(x_1,\ldots,x_j,\ldots,x_n)\leq f(z_1,\ldots,z_j,\ldots,z_n)\leq f(y_1,\ldots,y_j,\ldots,y_n).
\end{align}
finally, we have to for all $u=(u_1,\ldots,u_n)$ in $R$,
\begin{align}
  f\left( \prod_{j=1}^nx_j\right)\leq f\left( \prod_{j=1}^nu_j\right) \leq  f\left( \prod_{j=1}^ny_j\right)
\end{align}
or equivalent 
\begin{align}
f(R)\subset \left[  f\left( \prod_{j=1}^nx_j\right),f\left( \prod_{j=1}^ny_j\right)\right] .
\end{align}
Therefore
\begin{align}
f(R)= \left[  f\left( \prod_{j=1}^nx_j\right),f\left( \prod_{j=1}^ny_j\right)\right] 
\end{align}
which means that  $\displaystyle f(R)=\varphi^{-1}\phi f\left(R\right) $.
 \begin{flushright}
$ \blacksquare$
 \end{flushright}
 An observation to the above theorem, is that it remains true if the gradient is annulled at the most in two  vertexes   so that gradient does not cancel out in other points of $R$. For example, there is no problem if it is only annulled in a single vertex , in the case of two vertex, one of the points must correspond to a local maximum and the other a local minimum, or also that one of the two points is a saddle point. Now it can not happen that there are two maximum points or two minimum points, because, this would imply that there must exist a point of R between those points such that the gradient is zero. To this condition we will call \textit{free of sigularity except for the most in two vertexes.}
  
\begin{corollary} Under the same hypothesis of the above theorem.
Let $R=\displaystyle\bigcup_{j=1}^mR_j$ where $R_j$ are free of sigularity except for the most in two vertexes, then
\begin{align}
    f(R)=\displaystyle\bigcup_{j=1}^m\varphi^{-1}\phi f\left(R_j\right).
\end{align}
 
\end{corollary}
\textbf{Proof} Indeed $f(R)=f\left(\displaystyle\bigcup_{j=1}^mR_j\right)=\displaystyle\bigcup_{j=1}^mf(R_j)=\displaystyle\bigcup_{j=1}^m\varphi^{-1}\phi f\left(R_j\right)$
 \begin{flushright}
$ \blacksquare$
 \end{flushright}
 \newpage
 A consequence of the theorem is the following corollary.
 \begin{corollary} Let the following sets;
\begin{enumerate}
    \item $K_c(\mathbb{R})_0^+=\{X\in K_c(\mathbb{R}),x_1\geq 0\}$,
    \item $K_c(\mathbb{R})_0^-=\{X\in K_c(\mathbb{R}),x_2\leq 0\}$,
     \item $K_c(\mathbb{R})^+=\{X\in K_c(\mathbb{R}),x_1> 0\}$,
    \item $K_c(\mathbb{R})^-=\{X\in K_c(\mathbb{R}),x_2< 0\}$.
\end{enumerate} Then
 \begin{enumerate}
    \item $\varphi(X+Y)=\varphi(X)+\varphi(Y)$ for all $X,Y\in Kc(\mathbb{R}) $,
    \item $\varphi(kX)=k\varphi(X)$ for all $X\in Kc(\mathbb{R})$ and $k\geq 0$,
    \item $\varphi(kX)=k\varphi(\widehat{X})$ for all $X\in Kc(\mathbb{R})$ and $k\leq 0$,
    \item $\varphi(XY)=\varphi(X)\varphi(Y)$ for all $X,Y\in Kc(\mathbb{R})_0^+ $,
    \item $\varphi(XY)=\varphi(X)\varphi(\widehat{Y})$ for all $X\in Kc(\mathbb{R})_0^- $ and $Y\in Kc(\mathbb{R})_0^+ $,
    \item $\varphi(XY)=\varphi(\widehat{Y})\varphi(\widehat{Y})$ for all $X,Y\in Kc(\mathbb{R})_0^- $.
\end{enumerate}
 \end{corollary}
 \textbf{Proof} By simple inspection.
 \begin{flushright}
$\blacksquare$
\end{flushright}
We can observe that the $\varphi$ application preserves the addition and multiplication for intervals intervals with a single sign. besides that we can observe
  the first point is that due to 2) and 3) we have that $\varphi$ is nonlinear. And on the other hand, in general these properties are not combinable, for example by combining the properties of additivity and multiplication, we obtain the following inconsistency $ \varphi ([1,2] (1- [0,1]) = \left( \begin{array}{cc}
0 & 0 \\ 
0 &  2
\end{array}\right)$ and on the other hand $ \varphi ([1, 2] - [1,2] [0,1]) = \left( \begin{array}{cc}
-1 & 0 \\ 
0 &  2
\end{array}\right)$, from where we come to a contradiction. 
\section{Resolution of Interval Equations}\label{resol}

 Suppose we have an interval equation, for example a linear equation $AX + B = C$, where all the components are intervals, what should be the procedure to solve this equation ?, assuming that there is some solution. we could for example consider the equation $ax + b = c$, where the values of this equation are defined over their corresponding intervals, that is to say that $a \in A$, and clear $x$ of the equation and then determine the image of the square region, using the fundamental theorem, however, what we will obtain is a region that contains the solution of the equation.
 
 \begin{example}
 Let the interval linear equation $[1,2]X+[0,1]=[1,3]$ we can verify that the solution is $[1,1]$ however, we can consider the function $f(a,b,x)=ax+b$ with $a\in[1,2]$ and $b\in[0,1]$ and solve the following equation $f(a,b,x)=c$ with $c\in[1,3]$ in terms of $x$, when clearing $x$ we get the following function$ g(a,b,c)=\frac{c-b}{a}$, now calculating the image of $g$ over the square region $[1,2]\times[0,1]\times[1,3]$ is $[0,3]$ which does not correspond to the solution of the equation. 
 \end{example}
 We will then give a theorem that gives us the procedure to determine the solution of an interval equation, however, this solution does not always exist, since, the matrix we obtain as a solution to the matrix equation associated with the equation does not always satisfy the condition of have the first entry less than or equal to the last entry.\\\\
 \newpage
 Below we present the main theorem of this article
 \begin{theorem}[Fundamental Theorem of Interval Equations]\label{teofumeq} Let $f:\mathcal{X}\subset\mathbb{R}^n\to\mathbb{R}$  analytic function,  $\displaystyle\prod_{i=1}^nX_j\subset\mathcal{X}$  free of singularity with $X_j=[a_j,b_j]$, and $X_0\subset f\left( \mathcal{X}\right) $ an compact interval. Suppose it exists a function  $g:\displaystyle\prod_{i=2}^nX_j\to\mathbb{R}$ be such a function that for all $x_0\in X_0$ exists $(x_2,\ldots,x_n)\in \displaystyle\prod_{i=2}^nX_j$ such that $f\left(g\left(\displaystyle\prod_{i=2}^nx_j\right),\displaystyle\prod_{i=2}^nx_j\right)= x_0$, consider $\widehat{X_1}$ the  switch $X_1$ with respect to $f$. Then the equation $f\displaystyle\left(\prod_{i=1}^nX_j\right) =X_0$ has  solution in $X_1$ if and only if  $\varphi\widehat{X_1}=g\left(\displaystyle\prod_{i=2}^n\varphi\widehat{X}_j\right)$ defines an matrix with the first entry less than or equal to the last entry, i.e that $a_1\leq b_1$.  
\end{theorem}
\textbf{Proof}  Consider  the following interval equation in $X_1$:
\begin{align}\label{eq1}
     f\displaystyle\left(X_1,\prod_{i=2}^nX_j\right) =X_0.
\end{align}
As $\displaystyle\prod_{i=1}^nX_j\subset\mathcal{X}$ is free of singularity, then 
\begin{align}
   f\displaystyle\left(X_1,\prod_{i=2}^nX_j\right)=\varphi^{-1}f\displaystyle\left(\varphi \widehat{X_1},\prod_{i=2}^n\varphi \widehat{X}_j\right)
\end{align}
Then $\displaystyle\varphi^{-1}f\left(\varphi \widehat{X_1},\prod_{i=2}^n\varphi \widehat{X}_j\right)=X_0$ or the equivalent
\begin{align}\label{eq3}
  \displaystyle f\left(\varphi \widehat{X_1},\prod_{i=2}^n\varphi \widehat{X}_j\right)=\varphi X_0  
\end{align}
On the other hand we have by hypothesis, that there exists a function $g:\displaystyle\prod_{i=2}^nX_j\to\mathbb{R}$ be such a function that for all $x_0\in X_0$ exists $(x_2,\ldots,x_n)\in \displaystyle\prod_{i=2}^nX_j$ such that $f\left(g\left(\displaystyle\prod_{i=2}^nx_j\right),\displaystyle\prod_{i=2}^nx_j\right)= x_0$, that is, we can clear  $X_1$ of the matrix equation, then 
\begin{align}\label{eq2}
    \varphi \widehat{X}_1=g\left(\displaystyle\prod_{i=2}^n\varphi \widehat{X}_j\right)
\end{align}
Suppose there is a solution to the equation (\ref{eq1}), then (\ref{eq2}) define a matrix solution for the equation  (\ref{eq3})  where the first entry less than or equal to the last entry.\\
Let's suppose that $ \varphi \widehat{X}_1=g\left(\displaystyle\prod_{i=2}^n\varphi \widehat{X}_j\right)$ defined a matrix such that the first entry less than or equal to the last entry,  by (\ref{eq3}), we have
\begin{align}
  \displaystyle f\left(g\left(\displaystyle\prod_{i=2}^n\varphi \widehat{X}_j\right),\prod_{i=2}^n\varphi \widehat{X}_j\right)=\varphi X_0      
\end{align}
then
\begin{align}
 \varphi^{-1}\displaystyle f\left(g\left(\displaystyle\prod_{i=2}^n\varphi \widehat{X}_j\right),\prod_{i=2}^n\varphi \widehat{X}_j\right)= X_0      
\end{align}
by theorem \ref{teofun}, we have
\begin{align}
   X_0=\varphi^{-1}\displaystyle f\left(g\left(\displaystyle\prod_{i=2}^n\varphi \widehat{X}_j\right),\prod_{i=2}^n\varphi \widehat{X}_j\right)= \displaystyle f\left( \varphi^{-1}  \widehat{g}\left(\displaystyle\prod_{i=2}^n\varphi \widehat{X}_j\right),\prod_{i=2}^nX_j\right)      
\end{align}
as $ \varphi \widehat{X}_1=g\left(\displaystyle\prod_{i=2}^n\varphi \widehat{X}_j\right)$ defines a matrix with the first entry less than or equal to the last entry, we have to $ \varphi^{-1}  \widehat{g}\left(\displaystyle\prod_{i=2}^n\varphi \widehat{X}_j\right)\in K_c(\mathbb{R})$, then, we have that the latter corresponds to a solution of the equation (\ref{eq1}).
\begin{flushright}
$\blacksquare$
\end{flushright}
we can generalize the previous theorem for square regions with an arbitrary amount of singularities (that is, points with some null component of the gradient).
\begin{corollary}\label{coroeq}
 Let $f:\mathcal{X}\subset\mathbb{R}^n\to\mathbb{R}$  analytic function,  $\displaystyle\prod_{i=1}^nX_j\subset\mathcal{X}$    with $X_j=[a_j,b_j]$, and $X_0\subset f\left( \mathcal{X}\right) $ an compact interval. Suppose it exists a function  $g:\displaystyle\prod_{i=2}^nX_j\to\mathbb{R}$ be such a function that for all $x_0\in X_0$ exists $(x_2,\ldots,x_n)\in \displaystyle\prod_{i=2}^nX_j$ such that $f\left(g\left(\displaystyle\prod_{i=2}^nx_j\right),\displaystyle\prod_{i=2}^nx_j\right)= x_0$. Let $\Gamma$ the interval in the larger variable $x_1$ such that $\Gamma\times\displaystyle\prod_{i=2}^nX_j\subset\mathcal{X}$ and let $\{R_{\alpha}\}\subset\Gamma$ and $\{R_{\beta}\}\subset\displaystyle\prod_{i=2}^nX_j$ with  such that:
 \begin{enumerate}
     \item $\displaystyle\bigcup_{\alpha,\beta}R_{\alpha}\times R_{\beta}=\Gamma\times\displaystyle\prod_{i=2}^nX_j$,
     \item each $\displaystyle\left\{\frac{\partial f}{\partial x_j}\right\}$  has a constant sign not null in int$(R_{\alpha}\times R_{\beta})$,
     \item $R_{\alpha}\times R_{\beta}$ is free of sigularity except for the most in two vertexes.
 \end{enumerate}
 let's denote by $X_0^{\alpha,\beta}=f(R_{\alpha}\times R_{\beta})\cap X_0$ and $X_{\alpha}=X_1\cap R_{\alpha}$, let's take the following equation in $X_{\alpha}$
 \begin{align}\label{eqR}
     f(X_{\alpha}\times R_{\beta})=X_0^{\alpha,\beta}.
 \end{align}
 If $X_0\subset f\left(\Gamma\times\displaystyle\prod_{i=2}^nX_j\right)$, the following equation has a solution not empty
 \begin{align}\label{eqR2}
    f\left(X_{1}\times \displaystyle\prod_{i=2}^nX_j\right)=X_0 
 \end{align}
   if there is any solution of not empty for (\ref{eqR}) some $\alpha$ and $\beta$. Additionally, we have the solution of (\ref{eqR2}) if it exists, it is equal to $X_1=\displaystyle\bigcup_{\beta}\varphi^{-1}\phi\widehat{g}_{\alpha
  }\left(R_{\beta}\right)$ where $f(g_{\alpha}(R_{\beta}),R_{\beta})=X_0^{\alpha,\beta}$.
\end{corollary}
\textbf{Proof} Let's prove that $\displaystyle\bigcup_{\alpha}X_{\alpha}$ is a solution to the equation (\ref{eqR2}), for simplicity, we will say that in the case that the equation (\ref{eqR2}) some $\alpha$ and 
$\beta$ has no solution then we will say that $X_{\alpha}=\emptyset$, under the hypothesis of the corollary, we have
\begin{align}
&f\left(\bigcup_{\alpha}X_{\alpha}\times\displaystyle\prod_{i=2}^nX_j\right)=f\left(\bigcup_{\alpha}X_{\alpha}\times \bigcup_{\beta}R_{\beta}\right)= f\left(\bigcup_{\alpha,\beta}  X_{\alpha}\times R_{\beta} \right)\\
&=\bigcup_{\alpha,\beta}f\left(  X_{\alpha}\times  R_{\alpha} \right)= \bigcup_{\alpha,\beta} X_0^{\alpha,\beta}=\bigcup_{\alpha,\beta}f\left(R_{\alpha}\times R_{\beta}\right)\cap X_0 =f\left(\bigcup_{\alpha,\beta}R_{\alpha}\times R_{\beta}\right)\cap X_0 =X_0
\end{align}
then $\displaystyle\bigcup_{\alpha}X_{\alpha}$ is solution of  $f\left(X_1,\displaystyle\prod_{i=2}^nX_j\right)=X_0$. On the other hand as $R_{\alpha}\times R_{\beta}$ is sigularity except for the most in two vertexes, then by theorem \ref{teofumeq}, we have $X_{\alpha}=\varphi^{-1}\phi\widehat{g}_{\alpha}\left(R_{\beta}\right)$. Therefore $X_1=\displaystyle\bigcup_{\beta}\varphi^{-1}\phi\widehat{g}_{\alpha}\left(R_{\beta}\right)$.
\begin{flushright}
$\blacksquare$
\end{flushright}
\begin{example}
Consider the equation $[2,3]X+[7,10]=[-8,4]$. To solve this equation, let's consider the function $f:\mathbb{R}\times[2,3]\times[7,10]\to\mathbb{R}$ given by $f(x,a,b)=ax+b$, we have the following partial derivatives
\begin{align}
    \frac{\partial f}{\partial x}=a>0,   \frac{\partial f}{\partial a}=x \text{ and  }   \frac{\partial f}{\partial b}=1
\end{align}
then we have two regions free of singularity, in one the values of $x$ are  positivee and in the other they are negativ. in the positive case we have
\begin{align}
\left( \begin{array}{cc}
2 & 0 \\ 
0 &  3
\end{array}\right)\varphi X+\left( \begin{array}{cc}
7 & 0 \\ 
0 &  10
\end{array}\right) =
\left( \begin{array}{cc}
-8 & 0 \\ 
0 &  4
\end{array}\right)\Rightarrow X \in K_c(\mathbb{R})^-
\end{align}
contradiction. Therefore there is no positive solution. For the negative case we have
\begin{align}
\left( \begin{array}{cc}
3 & 0 \\ 
0 &  2
\end{array}\right)\varphi X+\left( \begin{array}{cc}
7 & 0 \\ 
0 &  10
\end{array}\right) =
\left( \begin{array}{cc}
-8 & 0 \\ 
0 &  4
\end{array}\right)\Rightarrow X=[-5,-3]
\end{align}
then for the corollary \ref{coroeq} we have that a solution for the general equation is $X=[-5,-3]$.
\end{example}
\begin{example}
Consider the  equation $X(1-X)=\displaystyle\left[0,\frac{1}{4}\right]$. To solve this equation, let's consider the function $f:\mathbb{R}\to\mathbb{R}$ given by $f(x)=x(1-x)$, the singularity-free regions of this function are $A=\displaystyle\left(-\infty,\frac{1}{2}\right)$ and $B=\displaystyle\left[\frac{1}{2},\infty\right)$ where we have that in $B$ the variable $X$ undergoes a switch. The solution of the restricted  equation in A is $\displaystyle\left[0,\frac{1}{2}\right)$ and the solution of the restricted  equation in $B$ is $\displaystyle\left[\frac{1}{2},1\right)$, by the Corollary \ref{coroeq} guarantees that a solution of the general equation is $\displaystyle\left[0,1\right]$. that we can observe that in particular the solution of the restricted equation to $B$ is also a solution of the general equation and the solution of the restriction $A$ adding the point $\displaystyle\frac{1}{2}$ is also a solution of the general equation.
\end{example}
\section{Another Approach}
There is another approach the resolution of interval equations apart from using square matrices, and it is to use the conscious set of polynomials of a variable over $(x ^ 2-x)$, it is easy to observe that this set is a ring with the operations of addition and multiplication of polynomials. The elements of this set are of the form $a + bh$, where $h ^ 2 = h$, which we will call pseudo complexes because of their similarity to complex numbers. In this case the application between the sets of the intervlar numbers and the complex pseudo numbers is given by $[a, b]\longmapsto a + (b-a) h$ and analogously as we did for the case of the square matrices, we can redo them for the set of the complex pseudo number.
\bibliographystyle{unsrt}  


\begin{thebibliography}{1}
 
\bibitem{moore}
 Ramon E. Moore.
\newblock  Method and applications of interval analysis.
\newblock  {\textit{Sism} }, (1979)

\

 
\end{thebibliography}

\end{document}